\documentclass[24pt]{amsart}
\usepackage{amsmath, amsthm, amscd, amsfonts}

\newtheorem{thm}{Theorem}[section]

\theoremstyle{definition}

\theoremstyle{remark}

\numberwithin{equation}{section}
\newfont{\kh}{msbm10}

\def\E{{\mathfrak E}}
\def\A{{\mathcal A}}

\begin{document}

\title[Sum of the reciprocal of the imaginary parts of the zeta
zeros]{Explicit approximation of the sum of the reciprocal of the
imaginary parts of the zeta zeros}

\author{Soheila Emamyari}

\address{Soheila Emamyari, \newline Department of Physics,
Institute for Advanced Studies in Basic Sciences, P.O. Box
45195-1159, Zanjan, Iran}

\email{emamyari@iasbs.ac.ir, soheila$_{-}$emamyari@yahoo.com}

\author{Mehdi Hassani}

\address{Mehdi Hassani, \newline Department of Mathematics, Institute for Advanced Studies in Basic
Sciences, P.O. Box 45195-1159, Zanjan, Iran}

\email{mhassani@iasbs.ac.ir, mmhassany@member.ams.org}



\begin{abstract}
In this note, we give some explicit upper and lower bounds for the
summation $\sum_{0<\gamma\leq T }\frac{1}{\gamma}$, where $\gamma$
is the imaginary part of nontrivial zeros $\rho=\beta+i\gamma$ of
$\zeta(s)$.
\end{abstract}

\maketitle


\section{Introduction}
The Riemann zeta function is defined  for $\Re(s)>1$ by
$\zeta(s)=\sum_{n=1}^{\infty}\frac{1}{n^s}$ and extended by analytic
continuation to the complex plan with one singularity at $s=1$; in
fact a simple pole with residues 1. The functional equation for this
function in symmetric form, is
$\pi^{-\frac{s}{2}}\Gamma\left(\frac{s}{2}\right)\zeta(s)=
\pi^{-\frac{1-s}{2}}\Gamma\left(\frac{1-s}{2}\right)\zeta(1-s)$,
where $\Gamma(s)=\int_0^{\infty}e^{-t}t^{s-1}dt$ is a meromorphic
function of the complex variable $s$, with simple poles at
$s=0,-1,-2,\cdots$ (see \cite{lebedev}). By this equation, trivial
zeros of $\zeta(s)$ are $s=-2,-4,-6,\cdots$. Also, it implies
symmetry of nontrivial zeros (other zeros $\rho=\beta+i\gamma$ which
have the property $0\leq\beta\leq 1$) according to the line
$\Re(s)=\frac{1}{2}$. The summation
$$
\A(T)=\sum_{0<\gamma\leq T }\frac{1}{\gamma},
$$
where $\gamma$ is the imaginary part of nontrivial zeros appears in
some explicit approximation of primes, and having some explicit
approximations of it can be useful for careful computations. This is
a summation over imaginary part of zeta zeros, and for approximating
such summations we use Stieljes integral and integrating by parts;
let $N(T)$ be the number of zeros $\rho$ of $\zeta(s)$ with
$0<\Im(\rho)\leq T$ and $0\leq\Re(\rho)\leq 1$. Then, supposing
$1<U\leq V$ and $\Phi(t)\in C^1(U,V)$ to be non-negative, we have
\begin{eqnarray}\label{ugammav}
\sum_{U<\gamma\leq V}\Phi(\gamma)=\int_{U}^V\Phi(t)dN(t)=-\int_{U}^V
N(t)\Phi'(t)dt+N(V)\Phi(V)-N(U)\Phi(U).
\end{eqnarray}
About $N(T)$, Riemann \cite{riemann} guessed that
\begin{equation}\label{nt}
N(T)=\frac{T}{2\pi}\log\frac{T}{2\pi}-\frac{T}{2\pi}+O(\log T).
\end{equation}
This conjecture of Riemann proved by H. von Mangoldt more than 30
years later \cite{davenport, ivic-zeta}. An immediate corollary of
above approximate formula, which is known as Riemann-van Mangoldt
formula is $\A(T)=O(\log^2 T)$, which follows by partial summation
from Riemann-van Mangoldt formula \cite{ivic-zeta}. In 1941, Rosser
\cite{rosser-41} introduced the following approximation of $N(T)$:
\begin{equation}\label{nt-r}
|N(T)-F(T)|\leq R(T)\hspace{10mm}(T\geq 2),
\end{equation}
where
$$
F(T)=\frac{T}{2\pi}\log\frac{T}{2\pi}-\frac{T}{2\pi}+\frac{7}{8},
$$
and
$$
R(T)=\frac{137}{1000}\log T+\frac{433}{1000}\log\log
T+\frac{397}{250}.
$$
This approximation allows us to make some explicit approximation of
$\A(T)$.
\section{Approximation of $\A(T)$}

\subsection{Approximate Estimation of $\A(T)$} As we set
above, $\gamma_1$ is the imaginary part of first nontrivial zero of
the Riemann zeta function in the upper half plane and computations
\cite{odlyzko} give us $\gamma_1=14.13472514\cdots$. On using
(\ref{ugammav}) with $\Phi(\gamma)=\frac{1}{\gamma}$, $0<U<\gamma_1$
and $V=T$, we obtain
\begin{equation}\label{at-int}
\A(T)=\int_U^T\frac{dN(t)}{t}=\int_U^T\frac{N(t)}{t^2}dt+\frac{N(T)}{T}.
\end{equation}
Substituting $N(T)$ from (\ref{nt}), we obtain
$$
\A(T)=\frac{1}{2\pi}\int_U^T\frac{\log\big(\frac{t}{2\pi}\big)}{t}dt-
\frac{1}{2\pi}\int_U^T\frac{dt}{t}+\frac{1}{2\pi}\log\frac{T}{2\pi}-\frac{1}{2\pi}
+O\left(\int_U^T\frac{\log(t)}{t^2}dt\right)+O\Big(\frac{\log
T}{T}\Big).
$$
Computing integrals and error terms, and then letting
$U\rightarrow\gamma_1^-$, we get the following approximation
$$
\A(T)=\frac{1}{4\pi}\log^2T-\frac{\log(2\pi)}{2\pi}\log T+O(1).
$$

\subsection{Explicit Estimation of $\A(T)$} Considering (\ref{nt-r})
and using (\ref{at-int}) with $2\leq U<\gamma_1$, for every $T\geq
2$ implies
$$
\int_U^T\frac{F(t)}{t^2}dt-\int_U^T\frac{R(t)}{t^2}dt+\frac{F(T)-R(T)}{T}\leq\A(T)\leq
\int_U^T\frac{F(t)}{t^2}dt+\int_U^T\frac{R(t)}{t^2}dt+\frac{F(T)+R(T)}{T}.
$$
A simple calculation, yields
$$
\frac{F(t)}{t^2}=\frac{d}{dt}\left\{\frac{1}{4\pi}\log^2t-\frac{1+\log(2\pi)}{2\pi}\log
t+\frac{\log^2(2\pi)-2\log(2\pi)}{4\pi}-\frac{7}{8t}\right\},
$$
and setting $\E(t)=\int_1^{\infty}\frac{ds}{st^s}$, we also have
$$
\frac{R(t)}{t^2}=\frac{d}{dt}\left\{-\frac{433}{1000}\frac{\log\log
t}{t}-\frac{137}{1000}\frac{\log
t}{t}-\frac{69}{40t}-\frac{433}{1000}\E(t)\right\}.
$$
The integral of $\E(t)$ converges for $t>1$; in fact
$\E(t)\sim\frac{1}{t\log t}$ when $t\rightarrow\infty$. Using the
fact $\frac{d}{dt}\E(t)=-\frac{1}{t^2\log t}$, we get
$$
\frac{1}{t\log t}-\frac{1}{t\log^2 t}<\E(t)<\frac{1}{t\log
t}-\frac{31}{95 t\log^2 t}
$$
for $t\geq 2$. Therefore, after letting $U\rightarrow\gamma_1^-$, we
obtain the following explicit upper bound
$$
\A(T)<\frac{1}{4\pi}\log^2T-\frac{\log(2\pi)}{2\pi}\log
T+\mathfrak{c_{au}}-\frac{137\log^2T+433\log
T-433}{1000T\log^2T}\hspace{10mm}(T\geq 2),
$$
where $\mathfrak{c_{au}}=0.43596427\cdots<\frac{109}{250}$, and an
easy computation verifies $-\frac{137\log^2T+433\log
T-433}{1000T\log^2T}<0$ for $T\geq 2.222$. Thus, we obtain
$$
\A(T)<\frac{1}{4\pi}\log^2T-\frac{\log(2\pi)}{2\pi}\log
T+\frac{109}{250}
$$
for $T\geq 2.222$. Similarly, we get
\begin{eqnarray*}
\A(T)&>&\frac{1}{4\pi}\log^2T-\frac{\log(2\pi)}{2\pi}\log
T+\mathfrak{c_{al}}\\&+&\frac{274\log^3T+866(\log\log
T)\log^2T+3313\log^2T+433\log
T-433}{1000T\log^2T}\hspace{10mm}(T\geq 2),
\end{eqnarray*}
where $\mathfrak{c_{al}}=0.06058187\cdots>\frac{3}{50}$ and for
$T\geq 2$ the last term in the above inequality is positive. So, we
obtain
$$
\A(T)>\frac{1}{4\pi}\log^2T-\frac{\log(2\pi)}{2\pi}\log
T+\frac{3}{50}
$$
for $T\geq 2$. Therefore we have proved the following result:
\begin{thm} Letting $\A(T)=\sum_{0<\gamma\leq T }\frac{1}{\gamma}$ with $\gamma$ is
imaginary part of zeta zeros, we have
\begin{equation}\label{at-lub}
\frac{15}{250}<\A(T)-\Big\{\frac{1}{4\pi}\log^2T-\frac{\log(2\pi)}{2\pi}\log
T\Big\}<\frac{109}{250},
\end{equation}
where the left hand side holds for $T\geq 2$ and the right hand side
holds for $T\geq 2.222$.
\end{thm}


\end{document}